\documentclass{amsart}
\usepackage{amssymb,amsmath,amsthm,amsfonts,amscd,amstext}
\usepackage[mathscr]{eucal}
\usepackage{array}
\usepackage{graphicx}
\usepackage{cite}
\newtheorem{lem}{Lemma}[section]
\newtheorem{prop}[lem]{Proposition}
\newtheorem{thm}[lem]{Theorem}

\newtheorem{df}{Definition}[section]
\newtheorem{rem}[lem]{Remark}

\begin{document}
\title{A note on the pentagram map and tropical geometry}
\author[T.~Kato]{Tsuyoshi Kato}
\address{Department of Mathematics, Graduate School of Science, Kyoto
University, Sakyo-ku, Kyoto 606--8502, Japan}
\email{tkato@math.kyoto-u.ac.jp}


%

\maketitle

\section{Pentagram map}
Theory of the pentagram map consists of  a beautiful integrable system,
which was originally  initiated by R.Schwartz
 from the view point of    the dynamical system on the set of convex polygons in projective geometry.
 Given an $n$-gon, the result of the action called the {\em pentagram map}, is also another $n$-gon given by the convex hull 
of the intersection points of consecutive shortest diagonals.
 These polygons  will shrink to a point
after many compositions, on the other hand if one considers
 the induced action on the set of convex polygons modulo projective transformations,
 then many interesting phenomena appear.
It  was
 extensively developed by 
V.Ovsienko-R.Schwartz-S.Tabachnikov, who discovered that 
 the most canonical object 
 from the view point of
projective geometry
is   the  {\em twisted polygons}.

From the view point of integrable system, 
the moduli space of twisted $n$-gons modulo projective transformations
constitutes complete feature, which equips with 
 the canonical Poisson structure, 
Casimirs and the integrals with the Hamiltonian, which commute with the pentagram map.
Moreover on the open subset consisted  by the universally twisted $n$-gons, the level set of the Hamiltonian 
is compact so that Liouville-Arnold-Jost theory can be applied.

From the view point of scale transform,
the continuous limit of the pentagram map $T$ turns out to be  the {\em Boussinesq equation} in integerble systems.
In the opposite direction to the scaling limit, we use 
{\em tropical geometry} which  is a scale transform 
between dynamical systems and 
provides with a correspondence between automata and real rational dynamics.
It  allows us to study  two very different dynamical systems at the same time and
 to induce some uniform analytic estimates by a comparison method.
 It eliminates detailed activities on rational dynamics
and extracts framework of their structure in automata.

 In this note we study dynamical properties of  the moduli spaces of twisted polygons
 in tropical geometry, and   induce 
{\em quasi-recursiveness} of the pentagram map,
which  allows errors from periodicity but within uniformly bounded amount
that are independent of the choice of initial values. 

Corresponding to the pentagram map, the pentagram automaton 
gives the dynamics whose orbits lie on the lower dimensional polytopes 
which is a tropical variety. As a future plan, 
it would be of interests for us to study the pentagram map form the 
view point of 
the interplay between integrable systems and tropical algebraic geometry.

\subsection{Twisted polygons}
We refer to $ [6,7,8,9] $ on foundations on theory of the pentagram map.
Here we quickly review some of the structure.

Let ${\frak P}_n$ be all the set of convex $n$-gons on the plane,
and denote the action   
$T: {\frak P}_n \to {\frak P}_n$
given by the convex hull 
of the intersection points of consecutive shortest diagonals.

A transformation that maps lines to lines is called a {\em projective transformation}.
Any projective transformation is given by a regular $3$ by $3$ matrix in homogeneous coordinates.
Any projective transformation acts on ${\frak  P}_n$, which commutes with $T$.
Let us denote
${\frak C}_n ={\frak P}_n / \sim$,
where $\sim$ is given by projective transformations, and 
denote the induced action by:
$$T: {\frak C}_n \to {\frak C}_n.$$

\vspace{3mm}

A {\em twisted $n$-gon} is a map
$\phi: {\mathbb Z} \to {\mathbb R P}^2$
which satisfies the equalities:
$$\phi(k+n) = M \circ \phi(k) \qquad (k \in {\mathbb Z})$$
for some projective automorphism $M$ called {\em monodromy}.
$\phi(k)$ correspond to  vertices on polygons, and so $M$ should be the identity
 for an untwisted polygon.

 Let us  put $v_i =\phi(i)$.
Two twisted $n$-gons $\phi_1, \phi_2$ are equivalent,
 if there is another projective transformation $\Psi$
such that $\Psi \circ \phi_1 =\phi_2$ holds.
Two monodromies satisfy the relation $M_2 = \Psi M_1 \Psi^{-1}$.

The pentagram map acts on the set of twisted polygons.
\begin{df}
Let us denote by ${\bf P}_n$ 
the space of twisted $n$-gons modulo equivalence, and denote the action by the pentagram map on it:
$$T: {\bf P}_n \to {\bf P}_n.$$
\end{df}
The pentagram maps is not fully defined on ${\bf P}_n$, but 
surely it does generically.

\subsection{Canonical coordinates}
Let us recall the {\em cross ratio} which is a classical invariant in projective geometry.
$\log$ of the cross ratio can also be used to describe the hyperbolic distance on the upper half plane.
For $t_1, \dots, t_4 \in {\mathbb R}$, the cross ratio is given by:
$$[t_1,t_2,t_3,t_4] = \frac{(t_1-t_2)(t_3-t_4)}{(t_1-t_3)(t_2-t_4)}.$$

It is well known that it is invariant under $PGl_2({\mathbb R})$ action.

\begin{df}
The canonical coordinate is given by:
\begin{align*}
& z_i = [v_{i-2}, v_{i-1}, ((v_{i-2},v_{i-1}) \cap (v_i, v_{i+1})), ((v_{i-2},v_{i-1}) \cap (v_{i+1}, v_{i+2}))],  \\
& w_i = [((v_{i-2},v_{i-1}) \cap (v_{i+1}, v_{i+2})), ((v_{i-1},v_i) \cap (v_{i+1}, v_{i+2})),
v_{i+1}, v_{i+2}]
\end{align*}
where $(v,v')$ is the line through $v$ and $v'$.
\end{df}

The next theorem gives the basic coordinates on the moduli space
of the twisted polygons. Later we study dynamical property of the pentagram map
passing through the coordinate.

\begin{thm}[6]
(1) The assignment:
$$\Phi: {\bf P}_n \to {\mathbb R}^{2n}$$
given by 
$$\phi \to \{(z_1, \dots, z_n),(w_1, \dots,w_n)\}$$
gives a generic coordinate so that it gives the local diffeomorphism.

(2)
The Pentagram map is given by:
$$T(z_i, w_i) = (z_i \frac{1-z_{i-1}w_{i-1}}{1-z_{i+1}w_{i+1}}, w_{i+1}  \frac{1-z_{i+2}w_{i+2}}{1-z_iw_i} ).$$
\end{thm}

\section{Main theorem}
We us use the metric on ${\mathbb R}^n$ given by:
$$d((x_1 , \dots, x_n) , (x'_1, \dots, x'_n) )
 \equiv \max_{1 \leq i \leq n} \{ |x_i -x'_i|\}$$
 which is quasi-isometric  to the standard one.

\begin{df}
The  dynamical system:
$$\varphi: {\mathbb R}^{2n} \to {\mathbb R}^{2n}$$
given by the formula mod $n$:
$$\varphi
\begin{pmatrix}
 x_i \\
 y_i
\end{pmatrix}
 = 
 \begin{pmatrix}
 x_i+
 \max(0, x_{i-1}+y_{i-1})-\max(0, x_{i+1}+y_{i+1}) \\
  y_{i+1}  +\max(0, x_{i+2}+y_{i+2})-\max(0,x_i+y_i)
  \end{pmatrix}
  $$
  is called the pentagram automaton.
  \end{df}
 Later we will see that the pentagram map is a kind of the dynamical framework 
 for the pentagram map, which can be seen at infinity (remark $3.3$).
  
  \vspace{3mm}

  For $k \geq 1$, let us put all the periodic points of the pentagram automaton:
  $$Per_k = \{ (\bar{x}, \bar{y})  \in {\mathbb R}^{2n} : \varphi^k (\bar{x}, \bar{y}) = (\bar{x}, \bar{y})\} .$$
Then  
for $t >1$, we put:
$$\tilde{Per}_k = \cup_{t>1} \{ (t^{x_1}, \dots, t^{x_n}, t^{y_1}, \dots, t^{y_n}) : 
(x_1, \dots,x_n,y_1, \dots,y_n) \in Per_k\} \subset {\mathbb R}^{2n}_{>0}.$$

\begin{df}
The quasi-periodic points of the pentagram map is given by:
\begin{align*}
{\bf Per}_k = \{ (- t^{x_1},   \dots,- t^{x_n}, t^{y_1}, \dots, & t^{y_n}) 
\cup  (  t^{x_1},  \dots,t^{x_n}, -t^{y_1}, \dots, -t^{y_n}) : \\
& (t^{x_1}, \dots, t^{x_n}, t^{y_1}, \dots, t^{y_n}) \in \tilde{Per}_k\}
\end{align*}
for all $t >1$.
 \end{df}

\begin{lem}
Let us  put
$\bar{B}_n = \{(  \bar{z},    \bar{w}) \in {\mathbb R}^{2n}_{>0} : \   0<   z_i w_j  \leq 1\}$.
Then we have the inclusion:
$${\bf D}_n = \{(  \bar{z},    \bar{w}) : 
(  \bar{z},  -  \bar{w}) \in \bar{B}_n  \} \cup  
 \{(  \bar{z},    \bar{w}) : 
( -  \bar{z},    \bar{w}) \in \bar{B}_n  \} \subset {\bf Per}_n.$$
\end{lem}
{\em Proof:}
The inclusion
$B_n = \{  (\bar{x}, \bar{y}) : x_i +y_j \leq 0\} \subset Per_n$
is satisfied, since
 $\varphi^l(x_i,y_i) = (x_i, y_{i+l})$ hold
for all $l$.

Then the condition $x_i +y_j \leq 0$ is equivalent to
$0< t^{x_i} t^{y_j} \leq 1$.

This completes the proof.

\vspace{3mm}

Our main result is the following:
\begin{thm}
Let $\phi$ be a twisted $n$-gon and 
$(\bar{z}_0, \bar{w}_0) \in {\bf Per}_k $ be any initial value.

Let us 
consider the orbit by the Pentagram map $T$ and denote:
$$T^m(\bar{z}_0, \bar{w}_0)=
(z_1^m, \dots,z_n^m, w_1^m, \dots, w_n^m) .$$
Then the uniform estimates hold for $1 \leq i \leq n$:
$$ 0 \ < \ |\frac{z_i}{z_i^k}|^{\pm 1} , \quad  |\frac{w_i}{w_i^k}|^{\pm 1}  \  \leq \  4^{(5^k-1)/4}.$$
\end{thm}

At infinity of the pentagram map, we can see the discrete integrable systems.
We have an application of integrability of the pentagram map passing through
tropical transform:

\begin{thm}
Any orbits of the  pentagram automaton
$\varphi: {\mathbb R}^{2n} \to {\mathbb R}^{2n}$
 lie generically on the $2[\frac{n+1}{2}] -2$ dimensional polytopes.
\end{thm}
The proofs of the above two theorems are   given in the next section.

\vspace{3mm}

\begin{rem}
There are particular class of twisted $n$-gons called 
universally convex $n$-gons.
Let ${\bf U}_n$ be the set of universally convex twisted $n$-gons modulo projective equivalence.
It  is known to be  open in ${\bf P}_n$ and
 is invariant under the action $T$.
Combination of  theorem $2.2$ with $2.4$ below suggest the `sizes' of the invariant tori,
which cannot be so `large' compared to the action of the pentagram map.
One might expect to obtain suitable notions concerning this.
\end{rem}

\begin{thm}[6]
Almost every point on ${\bf U}_n$ lies on a smooth torus that has  $T$-invariant affine structure.
\end{thm}
This follows  from complete integrability over ${\bf U}_n$
applied to Arnold-Liouville-Jost theorem.

\section{Tropical geometry}
\subsection{Tropical transform}
A relative $(\max, +)$-function $\varphi$ 
is a  piecewise linear function of the form: 
$$\varphi(\bar{x})=   
 \max(\alpha_1 + \bar{a}_1 \bar{x}, \dots , \alpha_m+ \bar{a}_m \bar{x})
 - \max(\beta_1 + \bar{b}_1 \bar{x}, \dots , \beta_l+ \bar{b}_l \bar{x})$$
 where
 $\bar{a}_k \bar{x}= \Sigma_{i=1}^n a_k^i x_i$, 
 $\bar{x} =(x_1, \dots,x_n) \in {\mathbb R}^n$,  
$ \bar{a}_k =(a_k^1, \dots, a_k^n) , \bar{b}_k \in {\mathbb Z}^n$
 and $\alpha_k , \beta_k  \in {\mathbb R}$. 
Notice that $\varphi$ is Lipschitz.

  We  say that the integer
$M \equiv ml$
 is the number of the components.

 For each relative $(\max, +)$ function $\varphi$ as above, we associate
a  parametrized rational function,
which we call the relative elementary function:
$$f_t(\bar{z}) = 
\frac{\Sigma_{k=1}^m t^{\alpha_k} \bar{z}^{\bar{a}_k}}
{\Sigma_{k=1}^l t^{\beta_k} \bar{z}^{\bar{b}_k}} \qquad (\bar{z}=(z_1, \dots, z_n) \in {\mathbb R}^n_{>0})$$
where
$  \bar{z}^{\bar{a}_k}= \Pi_{i=1}^n z_i^{a_k^i}$.
Notice that $f_t$ take positive values.

These two functions $\varphi$ and $f_t$ are connected passing through some intermediate
functions $\varphi_t$, which we 
 describe  shortly below.
For $t>1$, there is a family of semirings $R_t $ which are all 
the real number ${\mathbb R}$ as sets. 
The multiplications and the additions are respectively  given by:
$$x \oplus_t y = \log_t (t^x + t^y), \quad  x \otimes_t y = x+y.$$
As $t \to \infty$ one obtains  the equality
$x \oplus_{\infty} y = \max (x,y)$.

By use of $R_t$ as coefficients, one has 
{\em relative} $R_t$-{\em polynomials}: 
$$\varphi_t(\bar{x})=   
 (\alpha_1 + \bar{a}_1 \bar{x}) \oplus_t  \dots \oplus_t 
 ( \alpha_m+ \bar{a}_m \bar{x})
 - (\beta_1 + \bar{b}_1 \bar{x}) \oplus_t  \dots  \oplus_t ( \beta_l+ \bar{b}_l \bar{x})$$
The limit is given by the relative  $(\max ,+)$ function above:
$$\lim_{t \to \infty}  \varphi_t(\bar{x}) =
 \varphi(\bar{x}) .$$
Let us put $\text{Log}_t: {\mathbb R}^n_{>0} \to {\mathbb R}^n$ by
$(z_1, \dots, z_n) \to (\log_t z_1, \dots, \log_t z_n)$.
Then  $\varphi_t$ and  $f_t$ satisfy the following relation,
which can be verified by a straightforward calculation.
\begin{prop}[4,5,10]
$f_t \equiv (\log_t)^{-1} \circ \varphi_t \circ \text{Log}_t : {\mathbb R}_{>0}^n \to (0, \infty)$
is the relative elementary  function
$f_t(\bar{z}) = 
\Sigma_{k=1}^m t^{\alpha_k} \bar{z}^{\bar{a}_k}
/ \Sigma_{k=1}^l t^{\beta_k} \bar{z}^{\bar{b}_k}$.
\end{prop}

 These  three functions  $\varphi$, $\varphi_t$ and $f_t$
admit one to one correspondence between their presentations.
We  say that $\varphi$ is the corresponding $(\max,+)$-function to $f_t$.

Let us introduce the numbers, which satisfy the relation $cP_N(c) +1= P_{N+1}(c)$:
$$P_N(c) =  
\begin{cases}
 \frac{c^N-1}{c-1}  & c > 1, \\
  N  & c=1.
  \end{cases}$$

\subsection{Discrete dynamical systems}
Let us consider the discrete dynamical systems
$F: {\mathbb R}^{2n}_{>0} \to {\mathbb R}^{2n}_{>0} $
which is conjugate to the pentagram map, given by:
$$F(z_i, w_i) = (z_i \frac{1+z_{i-1}w_{i-1}}{1+z_{i+1}w_{i+1}}, w_{i+1}  \frac{1+z_{i+2}w_{i+2}}{1+z_iw_i} )
\quad (1 \leq i \leq n).$$
Let us consider its tropicalization with the  intermediate function:
$$\varphi,  \ \varphi_t: {\mathbb R}^{2n} \to {\mathbb R}^{2n}$$
given by:
\begin{align*}
& \varphi_t(\bar{x}, \bar{y}) = 
\begin{pmatrix}
&  x_i \otimes_t (0 \oplus_t  x_{i-1} \otimes_t y_{i-1}) - 0 \oplus_t  x_{i+1} \otimes_ty_{i+1} \\
& y_{i+1} \otimes_t  (0  \oplus_t   x_{i+2} \otimes_t y_{i+2}) - 0 \oplus_t  x_{i} \otimes_t y_{i}
\end{pmatrix}
\end{align*}
where $\varphi$ is the pentagram automaton.

Let us denote the orbits with the same initial value by:
\begin{align*}
& \varphi^l(x_1, \dots,x_n,y_1, \dots, y_n) = (x_{l,1}, \dots,x_{l,n},y_{l,1}, \dots,y_{l,n})  \\
& \varphi^l_t(x_1, \dots,x_n,y_1, \dots, y_n) = (x_{l,1}', \dots,x_{l,n}',y_{l,1}', \dots,y_{l,n}')  \\
\end{align*}

Let $(\bar{z}, \bar{w}) \in {\mathbb R}^{2n}_{>0}$ be an initial value, and consider the orbit:
$$F^l(\bar{z}, \bar{w})=(\bar{z}_l, \bar{w}_l).$$

\begin{lem}[1]
(1)
Let $M$ be the number of the components.
Then:
$$\sup_{\bar{q} \in {\mathbb R}^{2n}} |\varphi(\bar{q}) - \varphi_t(\bar{q} )| \leq \log_t M.$$

(2) Suppose the initial condition satisfies the relation:
$$x_i =\log_t z_i, \quad y_i = \log_t w_i$$
for all $1 \leq i \leq n$. Then the equalities hold:
$$x_{i,l}' = \log_t z_{i,l}, \quad y_{i,l}' = \log_t w_{i,l}$$
for all $l \geq 0$.
\end{lem}
{\em Proof:}
(2) follows from proposition $3.1$.

For convenience, we give a proof for (1).
We show the estimates:
$$|a_1 \oplus_t \dots \oplus_t a_m  - \max(a_1, \dots,a_m)| \leq \log_t m .$$
Assume  $a_1 =  \max(a_1, \dots,a_m)$. 
Then: 
\begin{align*}
a_1 \oplus_t  \dots \oplus_t a_m &  = \log_t (t^{a_1}+ \dots + t^{a_m})
= \log_t (t^{a_1}(1+ t^{a_2-a_1}  + \dots + t^{a_m-a_1} )) \\
& = a_1 + \log_t( 1+ t^{a_2-a_1} + \dots + t^{a_m-a_1})
\end{align*}
Since $a_i-a_1 \leq 0$ are  non positive, 
 the estimates hold:
 $$\log_t(1+ t^{a_2-a_1} + \dots + t^{a_m-a_1}) \leq \log_t m.$$

\begin{rem}
It follow from proposition $3.1$ and lemma $3.2$ that 
the pentagram automaton lies at infinity of  (conjugate of) the  pentagram map 
passing through the tropical transform.
\end{rem}

\begin{prop}
Let $c \geq 1$ and $M$ be the Lipschitz constant 
and the number of the components for $\varphi$ respectively.
Then the estimates hold:
$$|x_{i, l} -x_{i,l}'| , |y_{i, l} -y_{i,l}'|    \ \leq P_l(c)  \log_t M.$$
 \end{prop}
{\em Proof:}
For simplicity of the notation, let $q$ imply $x$ or $y$.

Firstly   one has the estimates by lemma $3.2$:
$$|q_{i,1} - q_{i,1}'  | \leq  \log_t M$$

Since $\varphi$ is $c$-Lipschitz,
the estimates hold:
\begin{align*}
&  |q_{i,2} - q_{i,2}'|  \leq | \varphi(\bar{q}_1) - \varphi_t(\bar{q}_1') |
 \leq | \varphi(\bar{q}_1) - \varphi(\bar{q}_1') | + 
| \varphi_t(\bar{q}'_1) - \varphi(\bar{q}_1') | \\
& \leq c  \ d(\bar{q}_1,  \bar{q}_1')+ \log_t M
\leq (c+1) \log_t M
\end{align*}
where the metric $d$ was given in section $2$.

Suppose the conclusion holds up to $l $. 
Then we have the estimates:
\begin{align*}
&  |q_{i,l+1} - q_{i,l+1}'|  \leq | \varphi(\bar{q}_l) - \varphi_t(\bar{q}_l') |
 \leq | \varphi(\bar{q}_l) - \varphi(\bar{q}_l') | + 
| \varphi_t(\bar{q}'_l) - \varphi(\bar{q}_l') | \\
& \leq c\ d(\bar{q}_l,  \bar{q}_l') + \log_t M
\leq (cP_l(c)+1) \log_t M =P_{l+1}(c) \log_t M.
\end{align*}

This completes the proof.

\subsection{Proof of theorem $2.2$}
 Let us choose 
  an initial value  $(\bar{z}, \bar{w}) \in \tilde{Per}_k$ 
and denote  the orbit
$F^l(\bar{z}, \bar{w})=(\bar{z}_l, \bar{w}_l)$.

Let us put the initial value:
$$x_i =x_i' = \log_t z_i, \quad y_i =y_i' =  \log_t w_i$$
for $1 \leq i \leq n$. 
By definition one obtains the periodicity:
$$\varphi^k(\bar{x}, \bar{y}) =( \bar{x}, \bar{y}) .$$

Now we have the estimates by proposition $3.3$:
$$d(\bar{q} , \bar{q}_k')  \leq d(\bar{q} , \bar{q}_k) + d(\bar{q}_k , \bar{q}_k')
= d(\bar{q}_k , \bar{q}_k') \leq \log_t M^{P_k(c)}.$$

Since  we have the equalities:
$$|x'_i - x_{i,k}' | = \log_t  (\frac{z_i}{z_i^k})^{\pm} , \quad 
|y'_i - y_{i,k}' | = \log_t  (\frac{w_i}{w_i^k})^{\pm} $$
and since $\log_t$ are monotone, we obtain the estimates:
$$(\frac{z_i}{z_i^k})^{\pm}, \quad (\frac{w_i}{w_i^k})^{\pm}  \leq M^{P_k(c)}.$$

Recall that the
 Pentagram map is given by:
$$T(z_i, w_i) = (z_i \frac{1-z_{i-1}w_{i-1}}{1-z_{i+1}w_{i+1}}, w_{i+1}  \frac{1-z_{i+2}w_{i+2}}{1-z_iw_i} ).$$
$T$ and $F$ can be transformed by change of the variables:
$$(\bar{z}, \bar{w}) \to (\bar{z}, - \bar{w}) \text{ or }  (- \bar{z},  \bar{w})$$

In the case of the pentagram map, the number of the component is $4$, and the Lipschitz constant is 
$3+2=5$.
This completes the proof.

\vspace{3mm}

\subsection{Proof of theorem $2.3$}
On the pentagram map, there are $2[\frac{n}{2}] +2$ algebraically independent invariants
$O_k, \quad E_k $ for 
$ 1 \leq k \leq [\frac{n}{2}]$ and $k =n$.

Let us consider the rational dynamics by $F$, which is given by exchange of the variable $z$ to $-z$.
Let us describe the concrete formulas for these invariants, 
and introduce the monomials:
$$Z_i = z_iw_iz_{i+1}$$

We say that 
$Z_i$ and $Z_j$ are consecutive, if $j \in \{i-2,i-1,i,i+1,i+2\}$,
$Z_i$ and $z_j$ are consecutive, if $j \in \{i-1,i,i+1,i+2\}$,
and $z_i$ and $z_{i+1}$ are consecutive.

An admissible monomial is given by:
$$O = Z_{i_1} \dots Z_{i_s} z_{j_1} \dots z_{j_r}$$
where no two factors are consecutive. Let us put:
$$|O| = s+t, \quad \text{ sign } (O) = (-1)^t$$

We have similar notions with $E$ by exchanging the role  between $z$ and $w$.

\begin{lem}[6]
The following list consists of all the conservation quantities:
\begin{align*} 
&  O_n= \Pi_{i=1}^n z_i ,  \quad  E_n= \Pi_{i=1}^n w_i , \\
& O_k = \Sigma_{|O|=k} O, \quad E_k = \Sigma_{|E|=k}  \text{ sign } (E)E.
 \end{align*}
\end{lem}

{\em Proof of theorem $2.3$:}
Let $\varphi, \varphi_t$ and $F$ be in $3.2$.
Let us choose any initial value $(x_1, \dots,x_l,y_1, \dots,y_l)$ and consider the orbits
$(x_{l.i},y_{l.i})$ and  $(x_{l.i}',y_{l.i}')$
for  $\varphi $ and $\varphi'$ 
respectively.

Let us put 
$z_i =t^{x_i}$ and $w_i = t^{y_i}$ for $1 \leq n$ and denote the orbit  as $(z^l_i,w^l_i)$ for $F$.
They satisfy the relations for all $l \geq 1$ and $1 \leq i \leq n$:
$$x_{l,i}' = \log_t z^l_i, \quad y_{l,i}' = \log_t w^l_i$$

There are constants $o_k >0$ and $e_k $ so that the invariants $O_k$ and $E_k$ 
take constant values under the action by $F$.

It follows from the above equalities that  the sums:
$$\Sigma_{i=1}^n \ x_{l,i} = \log_t c_n, \quad
\Sigma_{i=1}^n \ y_{l,i} = \log_t e_n$$
are both constants.

It follows from proposition $3.3$ that the estimates hold:

\begin{align*}
& |\Sigma_{i=1}^n \ x_{l,i} | \leq \Sigma_{i=1}^n \ |x_{l,i} -x_{l,i}'|  + |\Sigma_{i=1}^n \ x_{l,i}'| \leq
n P_l(c) \log_t M  +\log_t c_l \\
&  |\Sigma_{i=1}^n \ y_{l,i} | \leq \Sigma_{i=1}^n \ |y_{l,i} -y_{l,i}'|  + |\Sigma_{i=1}^n \ y_{l,i}'| \leq
n P_l(c) \log_t M  +\log_t e_l
\end{align*}
The left hand sides are both independent of $t >1$, and hence letting $t \to \infty$ on the right,
we obtain the equalities:
$$ \Sigma_{i=1}^n \ x_{l,i} = \Sigma_{i=1}^n \ y_{l,i} =0.$$
So we have induced two independent linear invariants.

Let us denote:
$$\tilde{O} =  X_{i_1} + \dots  + X_{i_s}  + x_{j_1} \dots  + x_{j_r}$$
where $X_i = x_i + y_i + x_{i+1}$, and we put $|\tilde{O}| \equiv |O|$ and equip with admissibility by the same way.
Then we take the maximum:
$$\tilde{O}_k = \max_{|\tilde{O}| =k} \{\tilde{O} \}$$
among all the admissible monomials of length $k$.

These are functions of the variables $(x_1, \dots,x_n, y_1, \dots, y_n)$.
It follows from lemma $3.2(1)$ with  the same argument as above, that  we have the invariants:
$$\tilde{O}_k(x_{l,1}, \dots, x_{l,n},y_{l,1}, \dots, y_{l,n}) =0$$

Let us consider the case of $E$. 
Because of the presence of the negative sign, we need to do some more.
Let us rewrite the invariants as:
$$E_k(1) \equiv \Sigma_{|E|=k, \text{ sign } E=1}   \ E = \Sigma_{|E|=k, \text{ sign } E =-1}  \ E \equiv E_k(-1) + e_k$$

We have similar expressions  $\tilde{E}_k(\pm 1)$.
Suppose $e_k >0$ holds. Then the orbit for $\varphi_t $ satisfy the equality:
$$E_k(1) (t^{x'_{l,1}}, \dots, t^{y'_{l,n}}) =  E_k(-1)  (t^{x'_{l,1}}, \dots, t^{y'_{l,n}})+ t^{ \log_t e_k}.$$
By letting $t \to \infty$, we obtain the equality:
$$\tilde{E}_k(1)(x_{l,1}, \dots, y_{l,n}) = \max( \tilde{E}_k(-1)(x_{l,1}, \dots, y_{l,n}) , 0).$$

Now we have the expression of the invariant:
$$\begin{cases}
\tilde{E}_k(1)(x_{l,1}, \dots, y_{l,n}) = \max( \tilde{E}_k(-1)(x_{l,1}, \dots, y_{l,n}) , 0) & e_k >0 \\
\tilde{E}_k(1)(x_{l,1}, \dots, y_{l,n}) =  \tilde{E}_k(-1)(x_{l,1}, \dots, y_{l,n})  & e_k =0 \\
\tilde{E}_k(-1)(x_{l,1}, \dots, y_{l,n}) = \max( \tilde{E}_k(1)(x_{l,1}, \dots, y_{l,n}) , 0) & e_k <0
\end{cases}$$

This completes the prof of theorem $2.3$.

\end{document}